\documentclass{article}
\usepackage{latexsym}
\usepackage{amsmath}
\usepackage{amsfonts}
\usepackage{amssymb}
\newcommand{\w}{\wedge}

  \newcommand{\be}{{\bf e}}
                           \newcommand{\e}{{\bf e}}

\newcommand{\cl}{C \kern -0.1em \ell}     %Clifford algebra

\newcommand{\ut}[1]{{\setbox0=\hbox{$#1$}\mathsurround=0pt
       \rlap{\raisebox{-0.8\dp0}{\raisebox{-0.8ex}
       {\kern -0.15ex\hbox{$\tiny\sim$}\kern 0.15ex}}}#1}}
\newcommand{\uti}[1]{{\setbox0=\hbox{$#1$}\mathsurround=0pt
       \rlap{\raisebox{-0.8\dp0}{\raisebox{-0.8ex}
       {\kern -0.3ex\hbox{$\tiny\sim$}\kern 0.3ex}}}#1}}

\newdimen\arrayruleHwidth                 % thick line
     \makeatletter
     \def\Hline{\noalign{\ifnum0=`}\fi\hrule \@height \arrayruleHwidth
         \futurelet \@tempa\@xhline}
     \makeatother

\input{psfig}

\newcommand{\MBS}{\raise 1pt \hbox{$\phantom{\scriptstyle>}$\space\space}}

\newcommand{\Y}[2]{Y^{#1}_{#2}} %macro command for Young operator
\newcommand{\G}[2]{G^{#1}_{#2}} %macro command for Garnir element
\newcommand{\Id}{{\bf 1}}
\newcommand{\mspn}{\mbox{\rm span}}
\newcommand{\rev}[1]{#1 \, \tilde{}}
\newcommand{\alphaq}{\alpha_q}

\newcommand{\ed}{\end{document}}

\newcommand{\reversion}[1]{#1 \, \tilde{}}

\newcommand{\BK}{\mathbb{K}}
\newcommand{\BF}{\mathbb{F}}
\newcommand{\BC}{\mathbb{C}}
\newcommand{\BR}{\mathbb{R}}
\newcommand{\BH}{\mathbb{H}}

\newcommand{\bff}{{\bf f}}

\newcommand{\fra}{\frac{1}{2}}

\newcommand{\beq}{\begin{equation}}
\newcommand{\eeq}{\end{equation}}

\begin{document}
\title{Hecke Algebra Representations in Ideals Generated by $q$-Young
Clifford Idempotents\thanks{Paper presented 
at the 5th International Conference on Clifford Algebras and their
Applications in Mathematical Physics, Ixtapa, Mexico, 
June 27 - July 4, 1999.}}
\author{Rafa{\l} Ab{\l}amowicz\\
Department of Mathematics Box 5054\\
Tennessee Technological University\\
Cookeville, TN 38505 USA\\
E-mail: rablamowicz@tntech.edu
\and
Bertfried Fauser\\
Universit\"at Konstanz\\
Fakult\"at f\"ur Physik, Fach M678\\
78457 Konstanz, Germany\\
E-mail: Bertfried.Fauser@uni-konstanz.de}
\date{August 12, 1999}
\maketitle

\begin{abstract}
It is a well known fact from the group theory that irreducible 
tensor representations of classical groups are suitably 
characterized by irreducible representations of the symmetric groups. 
However, due to their different nature, vector and spinor 
representations are only connected and not united in such description.

Clifford algebras are an ideal tool with which to describe symmetries 
of multi-particle systems since they contain spinor and vector 
representations within the same formalism, and, moreover, allow for a 
complete study of all classical Lie groups. 
%Unfortunately 
%there have been few attempts to introduce a symmetry classification 
%of such systems in the Clifford algebraic literature where mostly 
%single particle systems have been discussed. 
In this work, together 
with an accompanying work also presented at this conference, an 
analysis of $q$-symmetry -- for generic $q$'s -- based on the ordinary 
symmetric groups is given for the first time. We construct $q$-Young 
operators as Clifford idempotents and the Hecke algebra representations 
in ideals generated by these operators. Various relations as 
orthogonality of representations and completeness are given 
explicitly, and the symmetry types of representations is discussed. 
Appropriate $q$-Young diagrams and tableaux are given. The ordinary 
case of the symmetric group is obtained in  
the limit $q \rightarrow 1.$  All in all, a toolkit for Clifford 
algebraic treatment of multi-particle systems is provided. The 
distinguishing feature of this paper is that the Young operators 
of conjugated Young diagrams are related by Clifford reversion, 
connecting Clifford algebra and Hecke algebra features. This 
contrasts the purely Hecke algebraic approach of King and Wybourne, 
who do not embed Hecke algebras into Clifford algebras.

\noindent
{\bf MSCS: 15A66;  17B37;  20C30; 81R25 }

\end{abstract}

\section{Introduction}

\subsection{Motivation}

We investigate a possibility of implementing the symmetric group 
$S_n$ and its group algebra deformation, the Hecke algebra $H_\BF(n,q),$ 
as a subalgebra of the Clifford algebra of multivectors. The latter 
algebra is defined as the Clifford algebra of a bilinear form
with a suitably chosen anti-symmetric part. The presence of 
the antisymmetric part changes the structure of the corresponding 
Clifford algebra and allows one to introduce the needed deformation.

Our main interest in Clifford algebras arose from their ubiquitous 
appearance in mathematical physics, as it has been demonstrated 
many times by D.~Hestenes \cite{Hestenes-spacetime,Hestenes-foundmech}. 
Up to now, however, 
the main efforts have been devoted to the development of the real 
Dirac theory and other physical models such as, for example, the 
Weinberg--Salam theory \cite{Boudet} of electro-weak interactions. 
Despite this enormous range of applicability, there exist problems 
in mathematical physics not yet formulated or discussed in the 
Clifford algebra framework.

One major unsolved problem is the proper formulation of multi-particle 
theories. Quantum field theory is a theory of infinitely 
many particles which causes on one hand great problems with 
renormalization, but on the other it provides one of the most precise 
formalisms developed so far in physics. There have been only a few 
attempts to tackle the multiparticle problem \cite{Doran-states}, 
while Hestenes uses matrices in this case \cite{Hestenes-foundmech}. 
On the other hand, we succeeded in showing that quantum field theory, 
when treated in terms of generating functionals, can be reformulated 
by Clifford algebras 
\cite{Fauser-thesis,Fauser-positron,Fauser-transition}. However, to 
treat such complicated theories correctly, one is {\it forced\/} 
to introduce non-symmetric bilinear forms in Clifford algebras and  
there are at least two reasons why this needs to be done.

One reason is a problem of normal-ordering, which has to be performed 
in multiparticle quantum systems. The transition from time-ordered to 
normal-ordered generating functionals usually yields singularities, 
which can be seen as a calculational error if Clifford algebras are 
used \cite{Fauser-vertex}. The second reason is the connection of the 
vacuum structure of physical systems, which is intimately related 
to the antisymmetric part of the bilinear form in the Clifford 
algebra \cite{Fauser-vacua}. This is far beyond the abilities of 
other currently used methods.

The Clifford approach to such problems is very rigid. 
Treating Clifford numbers as single entities makes it easy to 
calculate with them but it hides the internal structure of the involved 
objects. Since essentially all physical observables
can be given in terms of two-spinors \cite{PenroseRindler}, we wish 
to have a mechanism which breaks up the Clifford numbers into smaller 
parts.

Taking an idea from the group theory, it is possible to characterize 
tensor products of irreducible representations of classical groups by 
the irreducible representations of the symmetric group, since both 
group actions commute. Such a situation, where one has a nontrivial 
action of $S_n,$ is then considered to be an $n$-particle system. 
The symmetric group plays hence a dominant role not only in mathematics, 
e.g., in combinatorics, but also in physics. As an example note that 
it was a group theoretical necessity to form mesons and hadrons 
from quarks after they had been postulated.

During the last two decades one has become aware that the deformed 
symmetric group algebra, i.e., the Hecke algebra, lays at the heart 
of the so-called quantized structures, see e.g. \cite{Connes,Majid}. 
In particular,  quantum groups have provided a powerful tool for 
solving 
lattice problems in statistical physics. Jones polynomial and 
the knot theory are related to Hecke algebras, see the extensive 
discussion in \cite{Fauser-hecke}. It is thus a quite natural idea 
to bring the symmetric group and its deformation, the Hecke algebra, 
into the Clifford formalism. Furthermore, the symmetric group is the 
Coxeter group of the Dynkin diagram of the ${\bf A}_n$ complex Lie
algebras \cite{Hiller} which explains the name given to the 
particular deformation discussed below.

To our knowledge, Clifford algebras have been used only marginally 
in the study of the symmetric groups. Only a spinor double cover 
is known, see \cite{HoffmanHumphreys}, however, D.~Finkelstein 
discussed this quite extensively in his lecture.

In this paper we show explicit calculations for 
$H_\BF(2,q)$ and $H_\BF(3,q)$ where $\BF$ is the base field 
of the group algebra in question. \footnote{%
We distinguish the field $\BF$ the algebras are built over and
the (double) field $\BK$, defined below, used in representations 
of Clifford algebras. They should not be confused.}

\subsection{Definitions}

Throughout this paper, we use Clifford algebras with an arbitrary not 
necessarily symmetric bilinear form $B.$ Such algebras can
be constructed using Chevalley's approach \cite{Chevalley}. However,
he utilizes only symmetric forms thereafter. This issue was clearly 
addressed in 
\cite{AblamowiczLounesto,Fauser-mandel,Oziewicz-FGTC,Oziewicz97},
while the connection to Hecke algebras was first discussed in 
\cite{Oziewicz95}. A mathematically sound approach to such
algebras including their applications can be found in 
\cite{Fauser-hecke,Fauser-transition}. There is a need imposed by 
quantum physics to use this type of Clifford algebras as it was shown 
in \cite{Fauser-vacua}.

Spinors are usually defined to be elements of a minimal left ideal 
of a Clifford algebra. Such ideals can be constructed as linear 
spaces generated by a primitive idempotent. This means, that the 
spinor space ${\cal S}=\,<\cl(B,V)\bff>_{\BK},$ where $\bff$ 
is a primitive idempotent and $\cl(B,V)$ is the full Clifford 
algebra of the pair $(B,V),$ where $B$ is the general bilinear form
and $V$ is a linear space over $\BF.$ It is well known, that the smallest 
faithful representation of a Clifford algebra is a spinor 
representation of dimension~$2^k,$ where $k$ is related to 
a Radon-Hurwitz number. Comparing the dimension of the space 
of endomorphisms of the spinor space ${\cal S}$ with that of the 
Clifford algebra, one finds several cases of representations over 
the field $\BK \cong \BR,\, \BC,\, \BH,\, \BR \oplus \BR $ or 
$\BH \oplus \BH.$ Moreover, since the primitive idempotents of 
a Clifford algebra decompose the unity 
$\Id = \bff_1 + \bff_2 + \cdots + \bff_k,$ 
one ends up with a $k$-dimensional right-linear space $\cal S$ 
over the appropriate field $\BK.$

The aim of this work is to provide a mechanism for breaking up the 
ordinary spinor representation of $\cl_{n,n}$ into tensor products 
of smaller representations using appropriate Young operators 
constructed as Clifford idempotents. This means, that a suitable 
Clifford algebra is used as a carrier space for various tensor product 
representations.

The Young operators for various Young diagrams provide us with a 
set of idempotents which decompose the unity $\Id$ of $\cl_{n,n}$ as
\begin{equation}
\Id = \Y{(\lambda_1)}{}+ \cdots + \Y{(\lambda_n)}{},
\end{equation}
where $(\lambda_i)$ is a partition of $n$ characterizing the 
appropriate Young tableau, that is, a Young diagram (frame) with 
an allowed numbering. We denote a Young tableaux by 
$\Y{(\lambda_i)}{i_1,\ldots,i_n},$ where $(\lambda_i)$ is an 
ordered partition of $n$ and $i_1,\ldots,i_n$ is an allowed 
numbering of the boxes in the Young diagram corresponding to 
$(\lambda_i)$ as in \cite{Hamermesh,Macdonald}. Furthermore, 
these Young operators are mutually annihilating idempotents
\begin{equation}
\Y{(\lambda_i)}{} \Y{(\lambda_j)}{} = 
\delta_{\lambda_i \lambda_j} \Y{(\lambda_j)}{}.
\end{equation}
It appears natural to ask if these Young operators can be used 
to give representations of the symmetric group {\em within\/} 
the Clifford algebraic framework. The representation spaces 
which appear as a natural outcome of the embedding of the symmetric 
group and its representations can then be looked at as multiparticle 
spinor states. However, these might not be spinors of the full 
Clifford algebra.

In order to be as general as possible, we give not only the 
representations of the symmetric group but also of the Hecke 
algebra $H_\BF(n,q).$ The Hecke algebra is the generalization 
of the group algebra of the symmetric group by adding the 
requirement that transpositions $t_i$ of {\em adjacent\/} elements 
$i,i+1$ are no longer involutions $s_i.$ We set $t_i^2 = (1-q) t_i + q$ 
which reduces to $s_i^2 = \Id$ in the limit $q\rightarrow 1.$

Hecke algebras are `truncated' braids, since a further 
relation (see (\ref{eq: t1}) below) is added to the braid group 
relations as in \cite{Artin}. A detailed treatment of this topic 
with important links to physics may be found, for example, in 
\cite{Goldschmidt,Wenzl} and in the references of \cite{Fauser-hecke}.

The defining relations of the Hecke algebra will be given according 
to Bourbaki \cite{Bourbaki}. Let $<\Id, t_1, \ldots, t_n >$ be a 
set of generators which fulfill these relations:
\begin{eqnarray}\label{eq: myt1}
t_i^2 &=& (1-q) t_i + q,                            \label{eq: t1} \\
t_i t_j &=& t_j t_i,   \quad \vert i-j \vert \ge 2, \label{eq: t2} \\
t_i t_{i+1} t_i &=& t_{i+1} t_i t_{i+1},            \label{eq: t3} 
\end{eqnarray}
then their algebraic span is the Hecke algebra. Since we will compare 
our results with those of King and Wybourne \cite{KingWybourne} 
-- hereafter denoted by KW -- we give a transformation to their 
generators $g_i,$ namely $g_i = -t_i,$ which results in a new quadratic 
relation
\begin{equation}
g_i^2 = (q-1) g_i +q
\end{equation}
while the other two remain unchanged. However, this small change in 
sign is responsible for great differences especially in the 
$q$-polynomials occuring in our formulas and in their formulas. 
One immediate 
consequence is that this transformation interchanges symmetrizers 
and antisymmetrizers. In particular, this replacement connects 
formula (3.4) in KW with our full symmetrizers while formula (3.3) in 
KW gives our full antisymmetrizers. Finally, the algebra morphism
$\rho$ which maps the Hecke algebra into the even part of an 
appropriate Clifford algebra can be found in \cite{Fauser-hecke}.

Let $<\Id, \be_1,\ldots, \be_{2n}>$ be a set of generators of 
the Clifford algebra $\cl(B,V),$ $V=\mspn\{\be_i\},$ with a 
non-symmetric $2n \times 2n$ bilinear form $B=[B(\be_i,\be_j)]=[B_{i,j}]$ 
defined as
\begin{equation}
B_{i,j} := \left\{
\begin{array}{cl}
  0,      & \mbox{if }\, 1\le i,j \le n           \,\mbox{ or }\,  n < i,j \le 2n, \\[1ex]
  q,      & \mbox{if }\, i = j-n                  \,\mbox{ or }\,  i-1-n = j,      \\[1ex]
  -(1+q), & \mbox{if }\, i+1=j-n                  \,\mbox{ or }\,  i=j+1-n,        \\[1ex] 
  -1,     & \mbox{if }\, \vert i-j-n \vert \ge 2  \,\mbox{ and }\, i>n,  \\[1ex]
  1,      & \mbox{otherwise. }
\end{array}
\right.
\label{eq:Bij}
\end{equation}
The most general case would have $\nu_{ij}\not=0$ in the last line of (\ref{eq:Bij}). 
For example, when $n=4,$ then
\begin{equation}
B = 
   \begin{pmatrix}
      0 & 0 & 0 & 0 & q &  -1-q & 1 & 1    \\
      0 & 0 & 0 & 0 &  -1-q & q & -1-q & 1 \\
      0 & 0 & 0 & 0 & 1 &  -1-q & q & -1-q \\
      0 & 0 & 0 & 0 & 1 & 1 & -1-q & q     \\
      1 & 1 & -1 & -1 & 0 & 0 & 0 & 0      \\
      q & 1 & 1 & -1 & 0 & 0 & 0 & 0       \\
     -1 & q & 1 & 1 & 0 & 0 & 0 & 0        \\
     -1 & -1 & q & 1 & 0 & 0 & 0 & 0
   \end{pmatrix}.
\label{eq:ourB}
\end{equation}
The bilinear form $B$ in (\ref{eq:ourB}) is our particular choice that guarantees 
that the following equations hold:
\begin{eqnarray}
\rho(t_i) &=& b_i := \be_i \wedge \be_{i+n},                        \label{eq:b1}\\
b_i b_j &=& b_j b_i, \; \text{ whenever } \; \vert i-j \vert \ge 2, \label{eq:b2}\\
b_i b_{i+1} b_i &=& b_{i+1} b_i b_{i+1}.                            \label{eq:b3}
\end{eqnarray}
This shows $\rho$ to be a homomorphism of algebras implementing the 
Hecke algebra structure in the Clifford algebra $\cl(B,V).$ One knows 
from \cite{Fauser-hecke} that $\rho$ is not injective, and that 
its kernel contains all Young diagrams which are not L-shaped 
(that is, diagrams with at most one row and/or one column).  The first 
instance, however, when this kernel is non-trivial occurs in $S_4$ where 
the partition $4=(2,2)$ gives a Young diagram of square form which is 
not L-shaped.

\section{The case of $H_{\BF}(2,q)$ and $S_2$}

We begin with $H_{\BF}(2,q)$ which reduces to $S_2$ in the limit 
$q\rightarrow 1.$ $H_{\BF}(2,q)$ is generated by $\{ \Id,b_1 \}.$ 
We have thus only one $q$-transposition, from which we can 
calculate a $q$-symmetrizer $R(12)$ and a $q$-antisymmetrizer $C(12).$

Notice that in the limit $q \rightarrow 1$ we have the following relations 
for a set of new generators defined as $s_i := \lim b_i \mbox{ when } q \rightarrow 1:$
\begin{enumerate}
\item [(i)]     $s_i^2 = 1,$ 
\item [(ii)]    $s_i s_j = s_j s_i,$ whenever $\vert i-j \vert \ge 2,$
\item [(iii)]   $s_i s_{i+1} s_i = s_{i+1} s_i s_{i+1},$
\item [(iii)']  $(s_i s_{i+1})^3 = \Id.$
\end{enumerate}
Property (iii)' follows from the fact that $s_i^2=1$ and $s_i^{-1} = s_i.$ 
This is a presentation of the symmetric group according to Coxeter-Moser
\cite{CoxeterMoser}. Now it is an easy matter to show that (ii) is valid 
for transpositions, and that (iii) can be calculated graphically using  
tangles as in Figure 1. 
%
%%%%%%%%%%%%%%%%%%%%%%%%%%%%%%%%%%%%%%%%%%%%%%%%%%%%%%%%%%%%%%
%Bertfried: if you don't have package seteps, you don't have %
%command \centereps available. Use then \centerline. I use   %
%\centereps because it shows the bounding box in .dvi file.  % 
%%%%%%%%%%%%%%%%%%%%%%%%%%%%%%%%%%%%%%%%%%%%%%%%%%%%%%%%%%%%%%
%
\begin{figure}[h]
\centerline{\psfig{figure=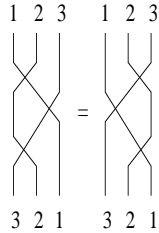,height=3.0cm,width=2.0cm}}
\caption{Tangles representing equation (\ref{eq:b3}).}
\label{fig:tangles}
\end{figure}
\noindent
In the cycle notation (iii) can be written as: 
$(12)(23)(12) = (23)(12)(23) = (13).$ 
In the Hecke case it does matter if one 'twists' the tangles one 
thread over or under the next one, which makes this relation quite 
nontrivial. It is sometimes called "the quantum-Yang-Baxter equation" 
in physics. In the case of the symmetric group, it does {\it not\/} 
matter which twist is used. In the crossings of the Hecke algebra we 
define the left-to-right moving tangle when reading from the top to 
the bottom as the {\it upper\/} one.

Thus, we define the $q$-symmetrizer $R(12)$ and the 
$q$-antisymmetrizer $C(12)$ (up to the normalization) as follows:
\begin{eqnarray}
R(12) & := & q + b_1 = q + \be_1 \wedge \be_5,     \label{eq:R12} \\
C(12) & := & \Id - b_1 = \Id - \be_1 \wedge \be_5. \label{eq:C12}
\end{eqnarray}
Notice, that the $q$-antisymmetrizer is related to the $q$-symmetrizer 
by the operation of reversion denoted by tilde 
in the Clifford algebra $\cl_{1,1},$ 
that is, $\rev{C(12)} = R(12)$ and $\rev{R(12)} = C(12).$ 
How do we know that $q + b_1$ gives the symmetrizer $R(12)$? Notice 
first that $R(12)$ is almost an idempotent since
\begin{equation}
R(12)R(12) = (1 + q) \be_1 \wedge \be_5 + q(1 + q ) = (1+q)R(12).
\end{equation}
Thus, when we normalize $R(12)$ by dividing it by $1+q,$ the new 
element denoted as $R(12)_q $ will be an idempotent.
$$
R(12)_{q}R(12)_{q} = \frac{(1 + q) \be_1 \wedge \be_5 + q(1 + q )}{(1 + q)^2} 
                   = \frac{\be_1 \wedge \be_5 +q}{1+q} 
                   = \frac{b_1 + q}{1+q} = R(12)_{q}.
$$
If we now take the limit of $R(12)_{q}$ as $q \rightarrow 1,$ we obtain
\begin{equation}
\lim_{q\rightarrow 1} R(12)_q = \frac{1+s_1}{2}
\end{equation}
with $s_1 = \be_1 \wedge \be_5$ squaring to $\Id$ (in the 
limit $q \rightarrow 1$) in agreement with (i) above. Then 
the expression $\frac{1}{2}(1+s_1)$ acts as a symmetrizer on, 
for example, functions of two variables. Likewise, the 
normalized $q$-antisymmetrizer gives an idempotent element 
\begin{equation}
C(12)_q = \rev{R(12)_q} =  \frac{1-b_1}{1+q}
\end{equation}
which in the limit $q \rightarrow 1$ gives the regular 
$S_2$ antisymmetrizer
\begin{equation}
\lim_{q\rightarrow 1} C(12)_q = \frac{1-s_1}{2}.
\end{equation}
Notice that $R(12)_q,\,C(12)_q$ and their limits 
$\frac{1}{2}(1+s_1),\, \frac{1}{2}(1-s_1)$ are pairwise 
mutually annihilating primitive idempotents in $\cl_{1,1}.$

Therefore, following the standard theory of Young operators, 
we define the $q$-Young operator $\Y{(2)}{1,2}$ as equal to the 
normalized $q$-symmetrizer $R(12)_q$ while the $q$-Young operator 
$\Y{(11)}{1,2}$ is defined as equal to the normalized 
$q$-antisymmetrizer $C(12)_q.$ We can conclude that $\Y{(2)}{1,2}$ 
and $\Y{(11)}{1,2}$ are mutually annihilating idempotents in 
$H_{\BF}(2,q) \subset \cl^{+}_{1,1} \subset \cl_{1,1}$ and that they 
decompose the unity
\begin{equation}
\Id = \Y{(2)}{1,2} + \Y{(11)}{1,2}. 
\end{equation} 

In the Clifford algebra $\cl_{1,1}$ we can construct various 
primitive idempotents. Since in our construction the Hecke 
algebra $H_{\BF}(2,q)$ is a subalgebra of the even part 
$\cl^+_{1,1}$ of $\cl_{1,1},$ the 
idempotents of the Hecke algebra must be even Clifford elements.  

It can be verified easily with CLIFFORD \cite{Ablamowicz,CLIFFORD} 
that the only two nontrivial mutually annihilating idempotents in 
the Hecke algebra $H_{\BF}(2,q)$ are the Young operators found 
above. In this case, the Young operators happen to be the two 
even primitive idempotents in the Clifford algebra $\cl_{1,1}.$
\footnote{The Clifford algebra $\cl_{1,1}$ is generated here by 
the $1$-vectors $\be_1$ and $\be_5.$ That is, we view $\cl_{1,1}$ 
as being embedded in the Clifford algebra $\cl_{4,4}$ of the 
$8 \times 8$ bilinear form $B$ given in (\ref{eq:ourB}).} 

With CLIFFORD we have also looked for any intertwining elements 
in the algebraic span of the two Young operators $\Y{(2)}{1,2}$ 
and $\Y{(11)}{1,2},$ and we have found none, as expected. That is, 
we have found no non-trivial element $T$ in the algebraic span 
of the Hecke algebra such that $T \Y{(2)}{1,2} = \Y{(11)}{1,2} T.$ 
Since the representations are only one-dimensional here, the 
Garnir element (see definition below) is zero.

\section{The case of $H_{\BF}(3,q)$ and $S_3$}

The Hecke algebra $H_{\BF}(3,q)$ is spanned by the basis elements 
$<\Id, b_1, b_2, b_{12}, b_{21}, b_{121}>$ which are expressed in 
terms of Grassmann polyomials as
\begin{equation}
\begin{array}{c}
b_1 = \be_1 \w \be_5, \quad b_2 = \be_2 \w \be_6, \\[0.75ex]
b_{12} = b_1 b_2 = -(1 + q) \Id + \be_1 \w \be_6 - \be_1 \w \be_2 \w \be_5 \w \be_6 + 
                    (1 + q) \e_2 \w \be_5, \\[0.75ex]
b_{21} = b_2 b_1 =-q(1+q)\Id + (1+q)\be_1 \w \be_6 - \be_1 \w \be_2 \w \be_5 \w \be_6 +
                   q\be_2 \w \be_5, \\[0.75ex]
b_{121} = b_1 b_2 b_1 = q\be_1\w\be_5-(1+2q)\Id+q\be_2\w\be_6+\be_1\w\be_6 \\[0.75ex]
                        \quad +(-1+q)\be_1\w\be_2\w\be_5\w\be_6-(-q-1+q^2)\be_2\w\be_5.\\
\end{array}
\end{equation}

We begin by defining our Young symmetrizer $\Y{(3)}{1,2,3}$ as in 
King and Wybourne:
\begin{equation}
\Y{(3)}{1,2,3} := \frac{q^3\Id + q^2 b_1 + q^2 b_2 + q b_{12} + q b_{21} + b_{121}}{(1+q+q^2)(1+q)}.
\end{equation}
By construction, our Young antisymmetrizer $\Y{(111)}{1,2,3}$ is 
defined as the reversion of the symmetrizer $\Y{(3)}{1,2,3},$ that is,
\begin{equation}
\Y{(111)}{1,2,3} := \reversion{\Y{(3)}{1,2,3}} =
\frac{\Id - b_1 - b_2 + b_{12} + b_{21} - b_{121}}{(1+q+q^2)(1+q)}.
\end{equation}
Also in this case King and Wybourne's full antisymmetrizer is the 
reversion of their full symmetrizer. However, the KW  Young operators
corresponding to Young tableaux which are conjugate of each other 
in the sense of Mcdonald \cite{Macdonald} (see also 
\cite{FultonHarris}) and generate representations of dimensions 
greater than one appear not to be related by the reversion. For 
example, King and Wybourne define
\begin{equation}
R(13) := q^3 \Id + b_{121} , \quad
C(13) := \Id - b_{121}
\end{equation}
yet $C(13) \neq \reversion{R(13)}$ since
\begin{equation}
\begin{array}{c}
C(13) - \reversion{R(13)}  =  2(q^2-q)\Id + (-1-q^2+2q)b_1 + (-1-q^2+2q)b_2 \\[1ex]
                             +(1-q)b_{12} + (1-q)b_{21},
\end{array}
\label{eq:C13R13}
\end{equation}
which does equal $0$ only in the limit $q \rightarrow 1,$ that is, 
not for an arbitrary $q.$ Therefore, the difference between the KW's 
antisymmetrizer $C(13)$ and their reversed symmetrizer 
$\reversion{R(13)}$ is the nonzero Hecke element (\ref{eq:C13R13}) which vanishes 
automatically in the limit 
$q\rightarrow 1.$ This is the reason for us to depart now from the 
KW formalism and introduce our own definitions of row-symmetrizers 
and column-antisymmetrizers, and require that they be related 
through the reversion. By doing so, we realize the conjugation of 
all Young tableaux as the Clifford algebra reversion. This is a 
crucial fact: since the reversion is the distinguished 
anti-automorphism of the Clifford algebra, we are able to connect 
it here to the automorphism of the Hecke algebra which interchanges 
the roots in the quadratic Hecke relation 
$(q,-1)\; \tilde \rightarrow \; (-1, q).$ Such interchange is an 
element of the Galois group \cite{Rotman}. We have thus established 
a direct connection between Clifford anti-involutions and the 
Galois group elements acting in the Hecke algebra. We consider our 
setting therefore to be in some sense natural.

To be as general as possible, we do the construction of the Young 
operators in two steps. First, we split the unit element $\Id$ in 
the Hecke algebra using the reversion into two non-primitive 
idempotents. Each of these idempotents generates a 
three-dimensional decomposable ideal. To achieve this goal, we can 
use any two of the following three equations since any two of them
imply the third:
\begin{eqnarray}
X+ X\tilde{~}  & = & \Id, \label{eq:X1}\\
X^2 & = & X,              \label{eq:X2}\\
X X\tilde{~} & = & 0      \label{eq:X3}.
\end{eqnarray}
By doing so, our goal is to find four Young operators known to 
exist from the general theory of the Hecke algebras for $n=3$ 
\cite{Goldschmidt,Wenzl}. The four $q$-Young operators will still 
have only one parameter and they will generalize four Young operators 
of $S_3$ described in Hamermesh \cite{Hamermesh} on p. 245. 
One of them will be a full symmetrizer, another one will be a full antisymmetrizer, 
and the other two will be of mixed symmetry. 

In the first step, we will find the most general element 
\begin{equation}
X = K_1\Id + K_2 b_1 + K_3 b_2 + K_4 b_{12} + K_5 b_{21} + K_6 b_{121}
\label{eq:genX}
\end{equation}
in the Hecke algebra $H_{\BF}(3,q)$ that satisfies~(\ref{eq:X1}). 
Upon substituting $X$ into~(\ref{eq:X1}) we have found that $X$ must have the 
following form:
\begin{equation}
\begin{array}{rcl}
X &=& (\fra q K_2 - \fra q K_6 + \fra + \fra q K_3 - \fra K_2 + \fra q^2 K_6 - \fra K_3) \Id  \\[1ex]
  & &  + K_2 b_1 + K_3 b_2 + K_4 b_{12} + (-K_4 + q K_6 - K_6) b_{21} + K_6 b_{121}.
\end{array}
\label{eq:Xsol1}
\end{equation}
The element $X$ in (\ref{eq:Xsol1}) belongs to a family parameterized by four real or 
complex parameters $K_2,K_3,K_4,K_6.$ Next we demand that $X$ also 
satisfies~(\ref{eq:X2}).

After substituting $X$ displayed in (\ref{eq:Xsol1}) into equation 
(\ref{eq:X2}), we have found six sets of solutions. The solutions are 
parameterized by complex numbers satisfying two similar but 
different quadratic equations:
\begin{equation} %alpha equation
\begin{array}{c}
(1+q)z^2+(-q^2 K_4+K_4+q K_2-1+K_2)z+K_4 K_2 + K_2^2-K_4-K_2\\[1ex]
+q K_4-q^2 K_4^2-q K_4^2-q^2 K_2 K_4+q K_2^2 = 0,
\end{array}
\label{eq:alphaeq}
\end{equation}
\begin{equation} %kappa equation
\begin{array}{c}
(1+q)z^2+(-q^2 K_4+K_4+q K_2+1+K_2)z+K_4 K_2+K_2^2+K_4+K_2\\[1ex]
-q K_4-q^2 K_4^2-q K_4^2-q^2 K_2 K_4+q K_2^2 = 0.
\end{array}
\label{eq:kappaeq}
\end{equation}
Let $\alpha$ be a root of equation (\ref{eq:alphaeq}) and 
$\kappa$ be a root of equation (\ref{eq:kappaeq}). We obtain the 
following six representatives $r_i,i=1,\ldots,6,$ of all solution 
families of (\ref{eq:X2}) which we again express in the Hecke 
algebra basis:\footnote{One might notice that the roots $\alpha$
and $\kappa$ provide us with two, possibly complex, solutions each, which 
yields $8$ real solutions: $4$ are linearly independent and $4$ are 
related by the reversion.}
\begin{eqnarray*}
r_1 & = & \frac{\Id}{1+q}-K_4 b_1+q K_4 b_2 + K_4 b_{12}- 
          \frac{(q^3 K_4+q+K_4 -1) b_{21}}{q(1+q)}       \\    
    &   & -\frac{(-K_4+q^2 K_4+1) b_{121}}{q(1+q)},         
\end{eqnarray*}
\begin{eqnarray*}
r_2 & = & \frac{q\Id}{1+q}-K_4 b_1+q K_4 b_2+K_4 b_{12}-
          \frac{(q^3 K_4-q+K_4+1) b_{21}}{q(1+q)}        \\
    &   & -\frac{(-K_4+q^2 K_4-1) b_{121}}{q(1+q)},          
\end{eqnarray*}
\begin{eqnarray*}
r_3 & = & \frac{\Id}{1+q}+q K_4 b_1-K_4 b_2 + K_4 b_{12}-
          \frac{(q^3 K_4+q+K_4-1) b_{21}}{q(1+q)}        \\
    &   & -\frac{(-K_4+q^2 K_4+1) b_{121}}{q(1+q)},         
\end{eqnarray*}
\begin{eqnarray*}
r_4 & = & \frac{q\Id}{1+q}+q K_4 b_1-K_4 b_2+K_4 b_{12}-
          \frac{(q^3 K_4-q+K_4+1) b_{21}}{q(1+q)}         \\
    &   & -\frac{(-K_4+q^2 K_4-1) b_{121}}{q(1+q)},         
\end{eqnarray*}
\begin{eqnarray*}
r_5 & = & \frac{\Id}{1+q}+K_2 b_1 + K_4 b_{12} \\
    &   & -\frac{(\kappa+K_2-q K_4+K_4-q^2 K_4^2+q K_2^2-q^2 K_2 K_4-
                q K_4^2+K_2^2+K_4 K_2) b_2}
               {(K_2+K_4-q K_4+\kappa)(1+q)} \\
    &   & -\frac{(q \kappa K_4+q K_2 K_4+q \kappa K_2-\kappa K_2) b_{21}}
               {q (K_2+K_4-q K_4+\kappa)} + 
          \frac{(\kappa K_2+q K_4^2) b_{121}}{q (-K_2-K_4+q K_4-\kappa)},
\end{eqnarray*}
\begin{eqnarray*}
r_6 & = & \frac{q \Id}{1+q}+K_2 b_1 + K_4 b_{12} \\
    &   & -\frac{(q^2 K_4^2-q K_2^2+q^2 K_2 K_4+q K_4^2+\alpha-q K_4+K_2+K_4-K_2^2-K_4 K_2)b_2}
               {(-K_2-K_4+q K_4-\alpha)(1+q)} \\
    &   & +\frac{(q K_2 K_4+q \alpha K_4+q\alpha K_2-\alpha K_2)b_{21}}
               {(-K_2-K_4+q K_4-\alpha)q} +    
          \frac{(\alpha K_2+q K_4^2)b_{121}}{q(-K_2-K_4+q K_4-\alpha)}.
\end{eqnarray*}
It can be checked with CLIFFORD that the rank of the set 
$\{r_i\},\,i=1,\ldots,6,$ is four. For our purpose we must select 
any four linearly independent elements, for example 
$\{r_1,r_2,r_3,r_5\},$ which we rename $\{f_1,f_2,f_3,f_4\}.$ 
It can also be verified with CLIFFORD that 
the elements $\{f_1,f_2,f_3,f_4\}$ satisfy the required relations 
(\ref{eq:X1}), (\ref{eq:X2}), and (\ref{eq:X3}).

We look for the Young operators obtained by one of the 
$f_i,i=1,\ldots,4.$ Due to the fact that the representation spaces 
which correspond to the symmetric $\Y{(3)}{1,2,3}$ and the 
antisymmetric $\Y{(111)}{1,2,3}$ Young operators respectively are 
one-dimensional, they cannot have any free parameters besides $q.$  
The full symmetrizer can be given according to KW as the $q$-weighted 
sum of all six Hecke basis elements. However, in our construction 
the full antisymmetrizer is defined as the reversion of the full 
symmetrizer, that is, $\Y{(111)}{1,2,3} := \reversion{\Y{(3)}{1,2,3}},$ 
as it was done in dimension two. Then we have:
\begin{eqnarray}
\Y{(3)}{1,2,3}   &:=& \frac{q^3\Id + q^2 b_1 + q^2 b_2 + q b_{12} + q b_{21} + b_{121}}
                           {(1+q+q^2)(1+q)}, \label{eq:Y3} \\[1ex]
\Y{(111)}{1,2,3} &:=& \frac{\Id - b_1 - b_2 + b_{12} + b_{21} - b_{121}}
                           {(1+q+q^2)(1+q)}. \label{eq:Y111}
\end{eqnarray}
Each of the four $f_i$ elements defined above generates a three-dimensional 
one-sided ideal in the Hecke algebra, or, in other words, together with 
its reversion $\tilde{f_i}$ it decomposes the unity in 
$H_{\BF}(3,q).$ Our objective is to further split each of these ideals 
into one- and two-dimensional vector spaces. The one-dimensional vector 
spaces will be generated by the full symmetrizer (\ref{eq:Y3}) and the 
full antisymmetrizer (\ref{eq:Y111}) respectively. The two-dimensional 
vector spaces will be generated by the mixed type Young operators. So 
each of the $f_i$ elements has to be a sum of a full (anti)symmetrizer 
and a Young operator of the mixed type. If we pick $f_1$ we notice that it 
must contain the full antisymmetrizer, because when the parameter $K_4$ 
is replaced with $1/(q+1)$ then the re-defined $f_1$ (or the $r_1$ 
defined above) reduces to an expression with alternating signs in the 
Hecke basis:
$$
\frac{\Id}{1+q}-\frac{b_1}{1+q}+\frac{q b_2}{1+q}+\frac{b_{12}}{1+q}-
\frac{q b_{21}}{1+q}-\frac{b_{121}}{1+q}.
$$
Therefore, by subtracting the full antisymmetrizer $\Y{(111)}{1,2,3}$ 
from $f_1$ we find our first Young operator $\Y{(21)}{1,3,2}$
of the mixed type: 
\begin{equation}
\begin{array}{rcl} 
\Y{(21)}{1,3,2} & = & f_1 - \Y{(111)}{1,2,3} \\[1.5ex]
              & = & \displaystyle{\frac{q \Id}{q+1+q^2}-
                    \frac{(q^3 K_4+2 q^2 K_4+2 q  K_4-1+ K_4) b1}{q^3+2 q^2+2 q+1}} \\[1.5ex]
              &   & +\displaystyle{\frac{(K_4 q^4+2 q^3 K_4+2 q^2 K_4+q K_4+1)b_2}
                                  {(q^3+2 q^2+2 q+1)}} \\[1.5ex]
              &   & +\displaystyle{\frac{(q^3 K_4+2 q^2 K_4+2 q K_4-1 + K_4) b_{12}}
                                  {(q^3+2 q^2+2 q+1)}} \\[1.5ex]
              &   & -\displaystyle{\frac{(K_4 q^5+ K_4 q^4+q^3 K_4+q^3+q^2 K_4+q K_4+q+ 
                                          K_4-1) b_{21}}
                                          {(q^3+2 q^2+2 q+1)q}} \\[1.5ex]
              &   & -\displaystyle{\frac{(K_4 q^4+q^3 K_4+q^2-q K_4+1-K_4) b_{121}}
                                       {(q+1+q^2) q (1+q)}}. \\
\end{array}
\label{eq:Y132}
\end{equation}
We define our second Young operator $\Y{(21)}{1,2,3}$ of the mixed type  as 
the reversion (conjugate) of $\Y{(21)}{1,3,2},$ that is, 
$\Y{(21)}{1,2,3} := \reversion{\Y{(21)}{1,3,2}},$ and we get:
\begin{equation}
\begin{array}{rcl} 
\Y{(21)}{1,2,3} 
&=& \displaystyle{\frac{q \Id}{q+1+q^2} +
                  \frac{(q^3 K_4-q^2+2 q^2 K_4+2 q K_4+ K_4) b_1}
                       {(1+q)(q+1+q^2)}} \\[1.5ex]
& & -\displaystyle{\frac{q(q^3 K_4+2 q^2 K_4+2 q K_4+q+ K_4) b_2}
                        {(1+q)(q+1+q^2)}} \\[1.5ex]
& & -\displaystyle{\frac{(q^3 K_4+2 q^2 K_4+2 q K_4+q+ K_4) b_{12}}
                        {q^3+2 q^2+2 q+1}} \\[1.5ex]
& & +\displaystyle{\frac{(K_4 q^5+ K_4 q^4+q^3 K_4+q^3-q^2+q^2 K_4+
                          q K_4-1+K_4)b_{21}}
                        {q (q^3+2 q^2+2 q+1)}} \\[1.5ex]
& & +\displaystyle{\frac{(K_4 q^4+q^3 K_4+q^2-q K_4+1-K_4)b_{121}}
                        {(q+1+q^2)q(1+q)}}. \\
\end{array}
\label{eq:Y123}
\end{equation}
Furthermore, the idempotent $\Y{(111)}{1,2,3}$ (resp. $\Y{(3)}{1,2,3}$) 
annihilates the idempotent $\Y{(21)}{1,3,2}$ (resp. $\Y{(21)}{1,2,3}$) 
when multiplied from both sides. This reflects the fact, that the 
representation spaces constructed in this way give a direct sum 
decomposition of the three-dimensional ideal generated by~$f_1$ 
(resp. $\tilde{f_1}$).  

Let's summarize the relationships between 
the Young operators:
\begin{equation}
\begin{array}{rcl}
f_1 = \Y{(111)}{1,2,3} +  \Y{(21)}{1,3,2} & , & \Y{(21)}{1,2,3} + \Y{(3)}{1,2,3} = \tilde{f_1},\\[1ex] 
\Y{(111)}{1,2,3}\Y{(21)}{1,3,2} = \Y{(21)}{1,3,2}\Y{(111)}{1,2,3} = 0 & , &
\Y{(21)}{1,2,3}\Y{(3)}{1,2,3} = \Y{(3)}{1,2,3}\Y{(21)}{1,2,3} = 0, \\[1ex]
\Y{(111)}{1,2,3} = \reversion{\Y{(3)}{1,2,3}} & , & \Y{(21)}{1,3,2} = \reversion{\Y{(21)}{1,2,3}}.  
\end{array}
\end{equation}
Thus, we have carefully built in the conjugation of the Young 
tableaux as the Clifford reversion. It can then be checked with 
CLIFFORD that the four Young operators displayed in (\ref{eq:Y3}), 
(\ref{eq:Y111}), (\ref{eq:Y132}), and (\ref{eq:Y123}) are primitive 
mutually annihilating idempotents which decompose the unity in the 
Hecke algebra since $f_1 + \tilde{f_1} = \Id.$ It can be easily 
verified that the Young operators of mixed type decompose into the 
row-symmetrizer and the colum-antisymmetrizer in accordance to 
Hamermesh \cite{Hamermesh} p. 245.  Our expressions however are 
different from those in KW. 

In order to represent our Young operator $\Y{(21)}{1,3,2}$ as a product 
of the row symmetrizer $R(13)$ and the column antisymmetrizer $C(12),$ 
we use previously defined $f_1 = r_1$ to define $C(12) := f_1$ and 
compute $R(13)$ from the equation
\begin{equation}
 Y^{(21)}_{1,3,2} = R(13) f_1 .
\label{eq:defR13}
\end{equation}
Notice that our $f_1=r_1$ is a generalization to $S_3$ of $C(12)$ 
from $S_2$ displayed in (\ref{eq:C12}). In an effort to be consistent 
with our previous discussion of $C(12)$ and $R(12),$ which were 
related by the reversion, we will later define 
$C(13) := \reversion{R(13)}$ and require that 
$R(13) + C(13) = R(13) + \reversion{R(13)} = \Id.$ 
Thus, when we solve (\ref{eq:defR13}) for $R(13),$ 
we get the following solution: 
\begin{equation}
\begin{array}{rcl}
R(13) & = & \displaystyle{\frac{q\Id}{1+q}}-
          \displaystyle{\frac{(-q^2+q^2 P_3+P_3 q-1+P_3) b_1}{(q+1+q^2)q}}+
          P_3 b_2 \\[1.5ex]
    &   & +\displaystyle{\frac{(q^2 P_3+P_3 q+P_3-1)b_{12}}{q(1+q+q^2)}} \\[1.5ex]
    &   & -\displaystyle{\frac{(q^5 P_3+q^4 P_3+q^3 P_3+q^2 P_3-q^2+P_3 q-1+P_3)b_{21}}
                   {(1+q) q^2 (1+q+q^2)}} \\[1.5ex]
    &   & -\displaystyle{\frac{(q^4 P_3+q^3 P_3-P_3 q+1-P_3)b_{121}}{(1+q) q^2 (1+q+q^2)}}
\end{array}
\end{equation}
where $P_3$ is an arbitrary real or complex parameter. It can be also 
verified with CLIFFORD that $R(13)$ is an idempotent element in the Hecke
algebra which is not a feature in KW. It is an interesting fact to note 
that the free parameter $P_3$ disappears when the product 
$R(13)f_1 = \Y{(21)}{1,3,2}$ is taken: the only free remaining parameter 
is $K_4.$ As a natural consequence of (\ref{eq:defR13}) we have
\begin{equation}
\Y{(21)}{1,2,3} = \tilde{f_1}\reversion{R(13)}.
\end{equation}
Obviously, we also have the following relations: 
\begin{equation}
C(13) := \reversion{R(13)}, \quad C(12):= f_1, 
\quad R(12) := \reversion{C(12)}.
\end{equation}

There are no non-trivially intertwining elements in the Hecke algebra 
which would connect the four Young operators except for the pair 
$\{\Y{(21)}{1,2,3}, \Y{(21)}{1,3,2}\}.$ That is, the only intertwiners 
that connect $\Y{(3)}{1,2,3}$ with $\Y{(111)}{1,2,3}$ actually 
annihilate them, and similarly for the pairs 
$\{\Y{(3)}{1,2,3}, \Y{(21)}{1,2,3}\}$ and  
$\{\Y{(111)}{1,2,3},\Y{(21)}{1,3,2}\}.$ On the other hand, there are 
many choices for an element $T$ in the Hecke algebra such that
\begin{equation}
T\Y{(21)}{1,2,3} = \Y{(21)}{1,3,2}T \neq 0.
\end{equation}
In fact, $T$ belongs to a five-parameter family of solutions. By 
assigning $0$ and $1$ values to four of those parameters, we have 
reduced the solutions to a one-parameter family parameterized only 
by~$K_4$ as follows: 
\begin{equation}
\begin{array}{rcl}
 T &=& \displaystyle{\frac{(1-K_4-q^3-K_4 q+K_4 q^3+q^4 K_4)\Id}
                    {q(2 K_4-q^2-q-q^3+2 q^4 K_4+8 K_4 q^2+6 K_4 q+6 K_4 q^3-1)}} \\[2ex]
   & & +\displaystyle{\frac{2(-1-q-q^2+K_4+2 K_4 q+K_4 q^3+2 K_4 q^2)b_1}
                    {q (2 K_4-q^2-q-q^3+2 q^4 K_4+8 K_4 q^2+6 K_4 q+6 K_4 q^3-1)}} \\[2ex]
   & & -\displaystyle{\frac{b_2}{1+q}} - \displaystyle{\frac{b_{12}}{q(1+q)}}+
       \displaystyle{\frac{b_{21}}{1+q}}+ \displaystyle{\frac{b_{121}}{q(1+q)}}
\end{array}
\end{equation}
where $K_4 \neq \displaystyle{\frac{(q^2+1)}{2(q^3+2 q^2+2 q+1)}}.$ 
Thus, the only intertwiners in the Hecke algebra which are not 
annihilators and which connect the mixed type Young operators constitute, 
in general, a four-parameter family of intertwiners while $T$ gives a 
one-parameter family.

We introduce now the Garnir elements $\G{(\lambda_i)}{i,j}$ (see KW) 
which will allow us to construct the six-dimensional representation 
spaces from the Young idempotents. Garnir elements can be seen to act 
as row or column permutations in Young tableaux, thereby generating 
non-standard tableaux which correspond to the basis vectors of 
the representation space. A {\it Garnir element\/} $\G{(21)}{1,1}$ 
has the following defining properties:
\begin{eqnarray}
\Y{(21)}{1,2,3}\,  \G{(21)}{1,1}   &   = & 0, \label{eq:G1} \\[1ex]
\G{(21)}{1,1}\,  \Y{(21)}{1,2,3}   & \neq& 0. \label{eq:G2} 
\end{eqnarray}
Our goal is to use a suitable Garnir element to decompose the 
three-dimensional one sided ideals generated by the Young operators 
of mixed symmetry $\Y{(21)}{1,3,2}$ and $\Y{(21)}{1,2,3}$ (see 
(\ref{eq:Y132}) and (\ref{eq:Y123}) respectively) into a direct sum 
of one-dimensional and two-dimensional vector spaces.  

When we first require that a general Hecke element $X$ shown in 
(\ref{eq:genX}) satisfies equation (\ref{eq:G1}) when substituted 
for $\G{(21)}{1,1},$ we find three different linearly independent 
solutions which we call $XX_1,\,XX_2,\, \mbox{and } XX_3.$ 
In the Hecke basis they look as follows:
\begin{equation}
\begin{array}{rcl}
XX_1 & = & (-K_6 q+K_2 q+K_4)\Id+K_2 b_1+ 
\displaystyle{\frac{t_1 b_2}{q}} + K_4 b_{12}+
             K_5 b_{21}+K_6 b_{121},
\end{array}
\label{eq:XX1}
\end{equation}
\begin{equation}
\begin{array}{rcl}
XX_2 & = & K_1 \Id+
          \displaystyle{\frac{t_2 b_1}{q^2 (1+q)}+ 
                        \frac{t_3 b_2}{q^3 (1+q)}-
                        \frac{b_{12}}{q(1+q)}}+
          K_5 b_{21}+K_6 b_{121},
\end{array}
\label{eq:XX2}
\end{equation}
\begin{equation}
\begin{array}{rcl}
XX_3 & = & K_1\Id+\displaystyle{
                  \frac{t_4 b_1}{q (q^3+2 q^2+2 q+1)}+
                  \frac{t_5 b_2}{q^2 (q^3+2 q^2+2 q+1)}} \\[2ex]
     &   &        -\displaystyle{\frac{(q-1) b_{12}}{q^3+2 q^2+2 q+1}}+
                   K_5 b_{21}+K_6 b_{121},
\end{array}
\label{eq:XX3}
\end{equation}
where the polynomial $t_1$ is parameterized by $K_2,K_4,K_5,K_6,$ 
and the polynomials $t_2,t_3,t_4,t_5$ are parameterized by 
$K_1,K_5,K_6.$\footnote{\mbox{~}
$t_1 = K_6 q^2+K_4 q^2-K_5 q-K_4 q-K_6 q+K_2+K_4,$ 
$t_2 = K_6 q^3+K_6 q^2+q^2 K_1+q K_1+1,$ 
$t_3 = q^5 K_6-q^4 K_5-q^3-q^3 K_5+q^2 K_1+K_6 q^2+q^2+q K_1-q+1,$
$t_4 = q-1+K_6 q+2 K_6 q^2+2 K_6 q^3+2 q K_1+2 q^2 K_1+K_1+q^3 K_1+q^4 K_6,$
$t_5 = q^6 K_6+q^5 K_6-q^5 K_5+q^4 K_6-q^4-2 q^4 K_5+
       K_6 q^3+q^3 K_1+2 q^3-2 q^3 K_5-K_5 q^2+K_6 q^2+
       2 q^2 K_1-2 q^2+2 q K_1+K_6 q+2 q+K_1-1.$}

With CLIFFORD it has been verified that the equation (\ref{eq:G2}) 
is satisfied automatically by each $XX_i,i=1,2,3,$ for all possible 
choices of the parameters. This means that we can use any of the three 
solutions as the Garnir element. Our choice is:
\begin{equation}
\begin{array}{rcl}
\G{(21)}{1,1} & = & (-K_6 q+K_2 q+K_4)\Id+K_2 b_1+ 
\displaystyle{\frac{t_1 b_2}{q}} + K_4 b_{12}+
             K_5 b_{21}+K_6 b_{121}.
\end{array}
\label{eq:G2111}
\end{equation}

Next we introduce the $q$-automorphism $\alphaq$ which replaces the 
reversion in changing the Garnir elements when acting from the 
left. This automorphism is in fact the inverse of the $b_i$'s 
generators and their products, the Hecke versors, and it is 
linearly extended to the entire Hecke algebra 
by means of the following definition:
\begin{equation}
\begin{array}{rcl}
\alphaq(b_{i_1} \ldots b_{i_s}) 
                      & = & 
        (\frac{-1}{q})^s (b_{i_1} \cdots b_{i_s})\tilde{~}       \\[1.5ex]
                      & = & 
        (\frac{-1}{q})^s (\tilde{b_{i_s}} \cdots \tilde{b_{i_1}})\\[1.5ex]
                      & = &  
\alphaq(b_{i_s}) \cdots \alphaq(b_{i_1}).
\end{array}
\label{eq:alphaq}
\end{equation}
For example,
$$
\begin{array}{c}
\alphaq(b_1) = \displaystyle{\frac{(q-1)\Id}{q} + \frac{b_1}{q}}, \quad
\alphaq(b_2) = \displaystyle{\frac{(q-1)\Id}{q} + \frac{b_2}{q}}, \\[2.0ex]
\alphaq(b_1 b_2) = \displaystyle{\frac{(1-2q+q^2)\Id}{q^2}+\frac{(q-1)b_1}{q^2}+
                                 \frac{(q-1)b_2}{q^2}+\frac{b_{21}}{q^2}}. \\[1.5ex]
\end{array}
$$
Note the fact that
\begin{equation}
\begin{array}{rcl}
\alphaq(b_i)b_i &=&  \frac{-1}{q} (\tilde{b_i} b_i) \\[1.5ex]
                &=&  \frac{-1}{q} (((1-q) - b_i) b_i) \\[1.5ex]
                &=&  \frac{-1}{q} ( (1-q) b_i -(1-q) b_i -q) \,=\, \Id.
\end{array}
\end{equation}
In fact, $\alphaq$ acts as the inverse on the generators and, as such, it 
depends on the presentation. However, while $\alphaq$ can be extended by the 
linearity to the entire Hecke algebra, it gives the inverse only of 
the products of the $b_i$'s (versors) and not of their sums. For example, 
\begin{equation}
\begin{array}{rcl}
\alphaq(\Id + b_1) = \displaystyle{\frac{(2q-1)\Id}{q} + \frac{b_1}{q}} 
& \neq &
(\Id + b_1)^{-1} = \displaystyle{\frac{(-2+q)\Id}{2(q-1)} + \frac{b_1}{2(q-1)}}, \, q \not= 1.
\end{array}
\end{equation}
One can check that it is not possible to use the reversion to 
connect the left and the right actions of Garnir elements even if 
this transformation connects the Young operators of conjugated 
diagrams. Looking at $\alphaq$ as an $q$-inverse, one could try 
to form a $q$-Clifford-Lipschitz group \cite{Fauser-q-groups} in 
the Hecke algebra by defining:
\begin{equation}
\Gamma_q := \{ X \,\vert\, \alphaq(X)X = \Id\,\}.
\end{equation}
When the automorphism $\alphaq$ is applied to the Garnir element 
$\G{(21)}{1,1},$ one gets
\begin{equation}
\begin{array}{rcl}
\alphaq(\G{(21)}{1,1}) & = & 
  -\displaystyle{\frac{t_6 \Id}{q^3} + \frac{t_7 b_1}{q^3}-\frac{t_8 b_2}{q^3}}+ 
   \displaystyle{\frac{(-K_6+K_6 q+K_5 q) b_{12}}{q^3}} \\[2.5ex]
                       &   & 
   +\displaystyle{\frac{(-K_6+K_4 q+K_6 q) b_{21}}{q^3}+\frac{K_6 b_{121}}{q^3}},
\end{array}
\end{equation}
where $t_6,t_7,t_8$ are polynomials.\footnote{\mbox{~}
$t_6 = K_2 q-2 K_6 q+K_6 q^3+K_6 q^2+K_5 q^2+
       K_6-q^3 K_2-q^4 K_2-K_5 q-q^4 K_4,$
$t_7 = K_6+q^2 K_2-K_5 q+K_5 q^2-2 K_6 q+K_4 q^2-K_4 q+K_6 q^2,$
$t_8 = -K_2 q+2 K_6 q-K_6+K_5 q-K_6 q^3-K_4 q^3.$}
With CLIFFORD we have verified that the six elements in the list 
$S$ below are linearly independent. As such, they provide a basis 
for the left regular representation of the Hecke algebra:
\begin{equation}
S = [\Y{(3)}{1,2,3},\,\Y{(21)}{1,2,3},\,\G{(21)}{1,1}\Y{(21)}{1,2,3},\,
  \alphaq(\G{(21)}{1,1})\Y{(21)}{1,3,2},\,\Y{(21)}{1,3,2},\,\Y{(111)}{1,2,3}].
\label{eq:Youngbasis}
\end{equation}
Recall that the original basis in the Hecke algebra was 
$[\Id,b_1,b_2,b_{12},b_{21},b_{121}].$ Each original basis 
element should be representable in terms of the new basis $S.$ 
Then we have:
$$
\Id = S_1+S_2+S_5+S_6,
$$
%S:=[Y3.123,Y21.123,&c(G21.1.1,Y21.123),&c(alpha2(G21.1.1),Y21.132),Y21.132,Y111.123]: 
%      S_1 ,  S_2  ,        S_3,                S_4,          S_5,      S_6
\begin{eqnarray*}
b_1 &=& S_1 +\frac{(K_4 q^3+2 K_4 q^2-q^2+2 K_4 q+K_4) S_2}{1+q} \\
    & & +\frac{p_1 q S_3}{p_2 (1+q)}+\frac{q^2 S_4}{p_3}+
        \frac{q(-K_4 q-K_2 q+K_6+K_5)S_5}{p_4}-q S_6,
\end{eqnarray*}
\begin{eqnarray*}
b_2 &=& S_1-\frac{q p_5 S_2}{1+q}-\frac{p_6 q S_3}{p_2 (1+q)}-
        \frac{q^3 S_4}{p_3}-\frac{p_7 S_5}{p_4}-q S_6,
\end{eqnarray*}
\begin{eqnarray*}
b_{12} &=& S_1 + \frac{p_8 S_2}{1+q}+ \frac{p_9 q S_3}{p_2 (1+q)} + 
           \frac{(1-q+q^2) q^2 S_4}{p_3}-\frac{q p_{10} S_5}{p_4}+q^2 S_6,
\end{eqnarray*}
\begin{eqnarray*}
b_{21} &=& S_1 -\frac{q p_5 S_2}{1+q}-\frac{p_{11} q^2 S_3}{p_2 (1+q)}-\frac{q^3 S_4}{p_3} - 
           \frac{q^2 (-K_4 q-K_2 q+K_6+K_5)S_5}{p_4}+q^2 S_6, 
\end{eqnarray*}
\begin{eqnarray*}
b_{121} &=& S_1 +\frac{p_{13} q S_2}{1+q}+\frac{p_{1,2} q^2 S_3}{p_2 (1+q)}+
            \frac{(q-1) q^3 S_4}{p_3}-\frac{q^2 p_{14} S_5}{p_4}-q^3 S_6 
\end{eqnarray*}
where to shorten the display we have introduced polynomials 
$p_i,\,i=1,\ldots,14$ (see the Appendix) and $S_i$ denotes the 
i-th element of the list $S,\,i=1,\ldots,6$. Finally, we can find, 
as expected, block-structured matrices of the basis elements $b_1$ 
and $b_2$ in the left regular representation (additional polynomials 
$p_i,\,i=15,\ldots,19,$ are also given in the Appendix). They are as follows:
\begin{equation}
M_{b_1} =  \begin{pmatrix}
1 & 0 & 0 & 0 & 0 & 0 \\
0 & {\displaystyle \frac {p_{15}}{1 + q}}  & 
    - {\displaystyle \frac {p_{16}}{q(1 + q)}}  & 0 & 0 & 0 \\ [2ex]
0 &   {\displaystyle \frac {q p_1}{p_{17}}}  
  &  - {\displaystyle \frac {p_{18}}{1 + q}}  & 0 & 0 & 0 \\ [2ex]
0 & 0 & 0 &  - {\displaystyle \frac {p_{19}}{p_4}}  
          &    {\displaystyle \frac {q^{2}}{p_3}} & 0 \\ [2ex]
0 & 0 & 0 &  - {\displaystyle \frac {p_{20}}{q\,p_4}}  
          &    {\displaystyle \frac {q p_{21}}{p_4}}  & 0 \\ [2ex]
0 & 0 & 0 & 0 & 0 &  - q
   \end{pmatrix}
\end{equation}
\begin{equation}
M_{b_2} =  \begin{pmatrix}
1 & 0 & 0 & 0 & 0 & 0 \\
0 &  - {\displaystyle \frac {q p_{5}}{1 + q}}  
  &    {\displaystyle \frac {p_{16}}{1 + q}}  
  & 0 & 0 & 0 \\[2ex]
0 & - {\displaystyle \frac {q p_6}{p_{17}}}  
  & {\displaystyle \frac {p_{21}}{1 + q}}  & 0 & 0
  & 0 \\ [2ex]
0 & 0 & 0 & {\displaystyle \frac {q p_{22}}{p_4}}  
          & {\displaystyle \frac {-q^{3}}{p_3}}  & 0 \\ [2ex]
0 & 0 & 0 &  - {\displaystyle \frac {p_{23}}{q p_4}}  
          &  - {\displaystyle \frac {p_7}{p_4}}  & 0 \\ [2ex]
0 & 0 & 0 & 0 & 0 &  - q
           \end{pmatrix}
\end{equation}
Matrices $M_{b_1}$ and $M_{b_2}$ satisfy, of course, the same quadratic 
Hecke relation (\ref{eq: t1}) as do $b_1$ and $b_2,$ and which happens 
to give the minimum polynomial $p(x) = x^2 - (1-q)x-q = (x-1)(x+q)$ 
for the basis elements and their matrix representations. The 
characteristic polynomial for the latter is $c(x) = (x-1)^{3}(x+q)^{3}.$ 
The trace of $M_{b_1}$ and $M_{b_2}$ is $3(q-1)$ while their 
determinants equal $-q^{3}.$ 

Finally, we build a left-regular matrix representation of the Young 
basis elements (\ref{eq:Youngbasis}) in the Young basis. 
\begin{equation}
\begin{array}{cc}
M_{S_1} =  \begin{pmatrix}
                       1 & 0 & 0 & 0 & 0 & 0 \\
                       0 & 0 & 0 & 0 & 0 & 0 \\
                       0 & 0 & 0 & 0 & 0 & 0 \\
                       0 & 0 & 0 & 0 & 0 & 0 \\
                       0 & 0 & 0 & 0 & 0 & 0 \\
                       0 & 0 & 0 & 0 & 0 & 0
           \end{pmatrix},
&
M_{S_2} = \begin{pmatrix}
      0 & 0 & 0 & 0 & 0 & 0 \\
      0 & 1 & 0 & 0 & 0 & 0 \\
      0 & 0 & 0 & 0 & 0 & 0 \\
      0 & 0 & 0 & 1 & 0 & 0 \\
      0 & 0 & 0 &  {\displaystyle\frac {p_{24}}{q^2}} & 0 & 0\\
      0 & 0 & 0 & 0 & 0 & 0
          \end{pmatrix},
\end{array}
\end{equation}
\begin{equation}
\begin{array}{cc}
M_{S_3} = \begin{pmatrix}
      0 & 0 & 0 & 0 & 0 & 0 \\
      0 & 0 & 0 & 0 & 0 & 0 \\
      0 & 1 & 0 & 0 & 0 & 0 \\
      0 & 0 & 0 & 0 & 0 & 0 \\
      0 & 0 & 0 &  - {\displaystyle \frac {p_{25}}{q^3}}  
                & 0 & 0 \\ [2ex]
      0 & 0 & 0 & 0 & 0 & 0
          \end{pmatrix},
&
M_{S_4} = \begin{pmatrix}
      0 & 0 & 0 & 0 & 0 & 0 \\
      0 & 0 &  - {\displaystyle \frac {p_{25}}{q^3}}  
                & 0 & 0 & 0 \\ [2ex]
      0 & 0 &  - {\displaystyle \frac {p_{24}}{q^2}}  
                & 0 & 0 & 0 \\ [2ex]
      0 & 0 & 0 & - {\displaystyle \frac {p_{24}}{q^2}}  
                    & 1 & 0 \\ [2ex]
      0 & 0 & 0 & 0 & 0 & 0 \\
      0 & 0 & 0 & 0 & 0 & 0
          \end{pmatrix},
\end{array}
\end{equation}
\begin{equation}
\begin{array}{cc}
M_{S_5} = \begin{pmatrix}
      0 & 0 & 0 & 0 & 0 & 0 \\
      0 & 0 & 0 & 0 & 0 & 0 \\
      0 & 0 & 1 & 0 & 0 & 0 \\
      0 & 0 & 0 & 0 & 0 & 0 \\
      0 & 0 & 0 & -{\displaystyle \frac {p_{24}}{q^2}}  
                    & 1 & 0 \\ [2ex]
      0 & 0 & 0 & 0 & 0 & 0
          \end{pmatrix},
&
M_{S_6} = \begin{pmatrix}
      0 & 0 & 0 & 0 & 0 & 0 \\
      0 & 0 & 0 & 0 & 0 & 0 \\
      0 & 0 & 0 & 0 & 0 & 0 \\
      0 & 0 & 0 & 0 & 0 & 0 \\
      0 & 0 & 0 & 0 & 0 & 0 \\
      0 & 0 & 0 & 0 & 0 & 1
         \end{pmatrix}
\end{array}.
\end{equation}

\section{Conclusions}

Motivated by the desire to describe symmetries of $q$-multiparticle 
systems we have constructed $q$-symmetrizers and $q$-antisymmetrizers 
in the Hecke algebras $H_\BF(2,q)$ and $H_\BF(3,q),$ and have related 
them by the reversion. That is, the Young operators constructed in 
the paper corresponding to the Young tableaux conjugate to each 
other in the sense of Macdonald and which generate representation 
spaces of dimension greater than one, have been related through 
the reversion in the Clifford algebra $\cl_{4,4}.$ This feature 
is not present in King and Wybourne. 

In $H_\BF(2,q)$ we found that the symmetrizer $R(12)$ and its 
reverse, the antisymmetrizer $C(12)$ were primitive idempotents 
in the Clifford algebra. We found no nontrivial intertwiners 
linking these two idempotents.

In $H_\BF(3,q)$ we first found four mutually annihilating idempotents 
splitting the unity in the algebra: two without parameters and two 
parameterized ones. The first two were the Young symmetrizer 
$\Y{(3)}{1,2,3},$ defined as in King and Wybourne, and the Young 
antisymmetrizer $\Y{(111)}{1,3,2},$ defined in this paper as the 
reverse of $\Y{(3)}{1,2,3}.$  The two parameterized idempotents were 
the Young operators of mixed symmetry 
$\Y{(21)}{1,2,3}$ and $\Y{(21)}{1,3,2},$ and they were also 
related by the reversion. We were able to factor 
$\Y{(21)}{1,3,2}$ into the row symmetrizer $R(13)$ and the 
column antisymmetrizer $C(12),$ that is, we were able to find 
$R(13)$ as an idempotent element in the Clifford algebra, a 
feature not found in King and Wybourne. Furthermore, we related 
the mixed-type Young operators through a five-parameter 
family of intertwiners. 

We have found a Garnir element $\G{(21)}{1,1}$ (in fact, three 
distinct families of such elements) which allowed us to further 
split the representation space of $H_\BF(3,q)$ from 
$3 \oplus 3$ to $1 \oplus 2 \oplus 2 \oplus 1:$ the 
one-dimensional spaces being generated by the Young symmetrizer 
$\Y{(3)}{1,2,3}$ and the Young antisymmetrizer $\Y{(111)}{1,3,2}$ 
while the two-dimensional spaces being generated by 
$\{\Y{(21)}{1,2,3}, \G{(21)}{1,1}\Y{(21)}{1,2,3}\}$ and 
$\{\alphaq(\G{(21)}{1,1})\Y{(21)}{1,3,2},\Y{(21)}{1,3,2}\}$ 
respectively. We introduced a Hecke algebra automorphism 
$\alphaq$ which acted as the inverse when applied to the 
Hecke basis elements: when applied to the Garnir element $\alphaq$ 
allowed us to generate the representation space of dimension 
six. The $\alphaq$ automorphism is expected to be useful in 
defining a $q$-Clifford-Lipschitz group in the Hecke algebra 
\cite{Fauser-q-groups}. Finally, we computed the matrix 
representation of the Hecke generators $b_1$ and $b_2$ in the 
Young basis. 

In the next step to be considered elsewhere we intend to extend 
this approach to $H_\BF(4,q).$ The connection with the spinor 
representations of the appropriate Clifford algebras in all 
three cases needs to be explored. Our approach to the Young 
operators related to the Young tableaux of various symmetries 
as Clifford (and Hecke) idempotents has been motivated by 
the need to describe the $q$-symmetries of multiparticle 
states. We have been able to define and construct all needed 
notations for a classification of tensor spaces w.r.t. the
$q$-symmetry. As mentioned above, the next natural step is to 
compute multiparticle spinors, that is, spin-tensors, as they are 
used for mesons or hadrons in QFT. Such spinors are, however, not
connected with the spinors of the Clifford algebras
used in this paper, but instead they should be connected with spinors 
of some appropriate sub-Clifford algebras and their tensor product. 
Our aim to break up a 'container Clifford algebra' into smaller pieces
will thus be achieved. We found that the Clifford algebra framework 
has proven to be most useful for this task. 

\section*{Appendix}

Polynomials below have been introduced as abbreviation to improve readability of 
the formulas displayed in the main text:
\begin{eqnarray*}
p_{1} &=&  q^4 K_4^2+3 q^3 K_4^2-K_4 q^3+4 q^2 K_4^2-K_4 q^2+3 q K_4^2 \\
      & &  -K_4 q-q+K_4^2-K_4,\\     
p_{2} &=&  K_4 q^4 K_6+q^4 K_4^2+q^3 K_4 K_6-q^3 K_4 K_2+K_4 q^2 \\
      & &  -q^2 K_2 K_4+q^2 K_4 K_6+K_6 q^2+q^2 K_4 K_5-q^2 K_4^2 \\
      & &  -K_2 q-K_6 q+q K_4 K_5-q K_4^2+K_4 q K_6-K_5 q \\
      & &  -q K_2 K_4-K_4 q-K_4^2-K_4 K_2+K_4+K_2, \\
p_{3} &=&  -K_6 q^3-K_4 q^3+K_6 q^2+K_5 q^2+K_6 q+K_5 q-K_6-K_4, \\
p_{4} &=&  -K_6 q^2-K_4 q^2+2 K_6 q+K_5 q+K_4 q-K_6-K_4, \\    
p_{5} &=&  K_4 q^3+2 K_4 q^2+q+2 K_4 q+K_4, \\    
p_{6} &=&  q^5 K_4^2+3 q^4 K_4^2+4 q^3 K_4^2+K_4 q^3+3 q^2 K_4^2 \\
      & &  +K_4 q^2+q K_4^2+K_4 q+K_4-1, \\    
p_{7} &=&  -K_4 q^3-q^3 K_2+2 K_6 q^2+K_5 q^2+K_4 q^2-2 K_6 q-K_5 q \\
      & &  -K_4 q+K_6+K_4,\\    
p_{8} &=&  q^5 K_4+q^4 K_4+q^3+K_4 q^3-q^2+K_4 q^2+K_4 q+K_4-1, \\    
p_{9} &=&  q^6 K_4^2+2 q^5 K_4^2+2 q^4 K_4^2+2 q^4 K_4+2 q^3 K_4^2 \\
      & &  +K_4 q^3+2 q^2 K_4^2+q^2+2 q K_4^2-K_4 q-q+K_4^2-2 K_4+1, \\    
p_{10} &=& q^3 K_2-K_6 q^3+K_6 q^2-q^2 K_2+K_5 q+K_2 q-K_5-K_6, \\     
p_{11} &=& q^4 K_4^2+3 q^3 K_4^2+4 q^2 K_4^2+K_4 q^2+3 q K_4^2 \\
       & & +K_4^2-K_4+1, \\     
p_{12} &=&  q^5 K_4^2+2 q^4 K_4^2+q^3 K_4^2+2 K_4 q^3-q^2 K_4^2 \\
       & &  +2 K_4 q^2-2 q K_4^2+2 K_4 q+q-K_4^2+2 K_4-1, \\     
p_{13} &=&  q^4 K_4+K_4 q^3+q^2-K_4 q-K_4+1, \\     
p_{14} &=&  q^2 K_2-K_6 q^2+K_6 q-K_2 q-K_4+K_5,\\ 
p_{15} &=&  K_4 q^3+2 K_4 q^2-q^2+2 K_4 q+K_4, \\     
p_{16} &=&  3 K_4 q^4 K_6+3 q^3 K_4 K_6+K_4 q K_6+2 q^2 K_4 K_6 \\
       & &  +2 q^5 K_4 K_6+2 q^2 K_4 K_5+K_6 q^4-2 q^4 K_4 K_2 \\
       & &  +q K_4 K_5+2 q^3 K_4 K_5-q^3 K_2-K_5 q-K_6 q-K_5 q^2 \\
       & &  -2 q K_2 K_4+K_2+q^6 K_4 K_6+q^4 K_4 K_5+K_4 \\
       & &  -3 q^2 K_2 K_4+q^4 K_4+K_4 q^2-K_4^2-3 q^3 K_4 K_2 \\
       & &  -q^5 K_4 K_2-2 q^3 K_4^2-q^3 K_5+q^6 K_4^2+q^5 K_4^2 \\
       & &  -3 q^2 K_4^2-2 q K_4^2-K_4 K_2,\\      
p_{17} &=& 2 K_4 q^4 K_6+2 q^3 K_4 K_6+K_4 q K_6+2 q^2 K_4 K_6 \\
       & & +q^5 K_4 K_6+2 q^2 K_4 K_5-q^4 K_4 K_2+q K_4 K_5 \\
       & & +q^3 K_4 K_5-K_5 q-K_6 q-K_5 q^2-2 q K_2 K_4+K_2 \\
       & & +q^4 K_4^2+K_4-2 q^2 K_2 K_4+K_6 q^3-q^2 K_2+K_4 q^3 \\
       & & -K_4^2-2 q^3 K_4 K_2-q^3 K_4^2+q^5 K_4^2-2 q^2 K_4^2 \\
       & & -2 q K_4^2-K_4 K_2, \\
\end{eqnarray*}

\begin{eqnarray*}     
p_{18} &=& K_4 q^3+2 K_4 q^2+2 K_4 q-1+K_4,\\     
p_{19} &=& -K_4 q^3-K_6 q^3+K_4 q^2+3 K_6 q^2+K_5 q^2-q^2 K_2 \\
       & & -2 K_4 q-2 K_6 q+K_6+K_4,\\     
p_{20} &=& K_4 q^4 K_6-3 q^3 K_4 K_6+2 K_4 q K_6-2 q^2 K_4 K_6 \\
       & & -q^5 K_4 K_6-q^2 K_4 K_5-K_6 K_4+q^4 K_4 K_2+q K_4 K_5 \\
       & & -q^3 K_4 K_5-q K_2 K_4+q^4 K_4^2+q^2 K_2 K_4-K_4^2 \\
       & & +2 K_6^2 q-K_6^2 q^2+2 q^3 K_4 K_2+K_5 K_4+q^5 K_4 K_2 \\
       & & +K_6 q^5 K_2-K_5 q^4 K_2-q K_2 K_6+K_6 K_5+q^3 K_2^2 \\
       & & -K_5 q^2 K_2+K_6 K_5 q+K_5 q^4 K_6+K_5 q^3 K_6 \\
       & & +2 K_6^2 q^4-K_6^2 q^5-2 q^3 K_2 K_5-2 q^4 K_2 K_6 \\
       & & -2 q^3 K_2 K_6+q^4 K_2^2,\\     
p_{21} &=& -K_4 q-K_2 q+K_6+K_5,\\     
p_{22} &=& -q^2 K_2+K_6 q^2-K_6 q-K_4 q+K_6+K_4,\\     
p_{23} &=& K_4 q^4 K_6-q^3 K_4 K_6+q^5 K_4 K_6-q^2 K_4 K_5 \\
       & & -K_6 K_4-q^4 K_4 K_2-q K_4 K_5-q^3 K_4 K_5+q K_2 K_4 \\
       & & +q^2 K_2 K_4-K_6^2-K_6^2 q^3+K_6^2 q^2-K_5 K_4-q^5 K_4 K_2 \\
       & & +K_6 q^5 K_2+K_5 q^4 K_2+q K_2 K_6-K_6 K_5-K_5 q^2 K_2 \\
       & & +K_6 K_5 q+2 K_6 q^2 K_5+K_5^2 q^2-q^5 K_2^2+K_5^2 q \\
       & & -K_5 q^4 K_6-K_5 q^3 K_6-K_6^2 q^4+2 q^4 K_2 K_6-q^4 K_2^2 \\
       & & +q^2 K_4^2+q K_4^2,\\     
p_{24} &=& -q^4 K_4^2+K_6-q^2 K_2+2 K_6 q^2-q^5 K_4 K_6 \\
       & & +2 q^3 K_4 K_5+q K_4 K_5-q^5 K_4^2-K_6 K_4-K_6 q-K_4 q^3 \\
       & & +q^2 K_4 K_6+K_4+K_5 q^2-q^2 K_4^2-q^3 K_2+q^3 K_4 K_6 \\
       & & +2 q^2 K_4 K_5+q^4 K_4 K_5-K_4 q-q K_4^2-q^3 K_4^2-K_4^2, \\     
p_{25} &=& -3 K_4 q^4 K_6+3 q^3 K_4 K_6+3 K_4 q K_6-3 q^2 K_4 K_6 +3 q^5 K_4 K_6 \\
& & -K_6 K_4-q^4 K_4 K_2+2 q K_4 K_5+2 q^3 K_4 K_5 +q K_2 K_4-2 q^6 K_4 K_6\\
& & -K_5^2 q^4-K_6^2 q^6-q^3 K_5^2 -K_4^3 q^8+q^4 K_4^3+q^5 K_4^3-2 q^4 K_4^2+2 K_6 q^5 K_5\\
& & +q K_4^3 +2 q^2 K_4^3+q^3 K_4^3+K_4 q^4 K_6^2+K_4^2 q^7 K_5+K_4^2 q^7 K_2 -2 K_4^2 q^8 K_6\\
& & -q^2 K_2 K_4+3 q^5 K_5 K_4 K_6-2 q^5 K_2 K_4 K_5 -q^4 K_2 K_4 K_6+3 K_5 q^3 K_4 K_6\\
& & -3 K_2 q^3 K_4 K_5+2 K_5^2 q^4 K_4+K_6 q^7 K_4 K_5+q^7 K_2 K_4 K_6-K_6^2 q^8 K_4\\
& & -3 K_2 q^4 K_4 K_5-K_2 q^3 K_4 K_6+4 K_5 q^4 K_4 K_6+K_6^2 q^3 -K_4^2+2 q^5 K_4 K_5\\
& & +K_6^2 q-2 K_6^2 q^2+q^3 K_4 K_2+q^5 K_4 K_2 +K_5^2 q^5 K_4-K_4^2 q^6 K_6+2 q^3 K_4^2 K_2\\
& & +2 q^4 K_4^2 K_2+2 q^5 K_4^2 K_2+q^6 K_4^2 K_2-4 q^3 K_4^2 K_6-3 K_4^2 q^4 K_6-K_4^2q^5K_6\\
& & -3 q^4 K_4^2 K_5-q^5 K_4^2 K_5-q^2 K_4^2K_6+2q^2 K_2 K_4^2-3 q^2 K_4^2 K_5-4 q^3 K_4^2 K_5\\
& & +q K_2 K_4^2+K_4^3+q^6 K_5 K_4 K_6+K_6 q^5 K_2+q^3 K_4^2-K_5 q^4 K_2+2 q K_2 K_6 \\
& & +K_6 K_5 q-2 K_6 q^2 K_5-K_5^2 q^2-K_5 q^4 K_6-K_6^2 q^4 +2 K_6^2 q^5-2 q^4 K_2 K_6\\
& & +q^3 K_2 K_6+K_4 q^5 K_6^2-q^6 K_4^2+q^5 K_4^2-2 q^2 K_5 K_4 K_2+q^2 K_5 K_4 K_6 \\
& & +K_6 K_4 K_2+K_4^2 K_2+K_4^2 K_6-2 q K_4^2 K_5-q K_4^2 K_6-K_6 K_2-q K_6 K_4 K_5 \\
& & -q K_5 K_2 K_4+2 K_5^2 K_4 q^3+K_5 K_2 q-q^6 K_5 K_2 K_4 \\
       & & -K_6 K_2 q^2+K_5^2 K_4 q^2-K_6^2 K_4 q-2 q^2 K_4^2+q K_4^2-K_4 K_2.
\end{eqnarray*}

\ed #quick \end{document}